# A New Algorithm for General Cyclic Heptadiagonal Linear Systems Using Sherman–Morrison–Woodbury formula

A.A. KARAWIA[1]
Computer Science Unit, Deanship of Educational Services,
Qassim University, Buraidah 51452, Saudi Arabia.
kraoieh@qu.edu.sa

## ABSTRACT

In this paper, a new efficient computational algorithm is presented for solving cyclic heptadiagonal linear systems based on using of heptadiagonal linear solver and Sherman–Morrison–Woodbury formula. The implementation of the algorithm using computer algebra systems (CAS) such as MAPLE and MATLAB is straightforward. Numerical example is presented for the sake of illustration.

## 1. INTRODUCTION

The cyclic heptadiagonal linear systems often occur in several fields such as numerical solution of ordinary and partial differential equations, interpolation problems, boundary value problems, etc. [1,2]. In many of these areas, it is necessary to obtain the solution of the cyclic heptadiagonal linear systems.

In this paper, we consider general cyclic heptadiagonal linear systems of the form

$$HX = R, \tag{1.1}$$

where

$$H = \begin{bmatrix} d_1 & a_1 & A_1 & C_1 & 0 & 0 & \cdots & & B_1 & b_1 \\ b_2 & d_2 & a_2 & A_2 & C_2 & \ddots & \ddots & \cdots & & B_2 \\ B_3 & b_3 & d_3 & a_3 & A_3 & \ddots & \ddots & \ddots & \cdots & \\ D_4 & B_4 & b_4 & d_4 & a_4 & \ddots & \ddots & \ddots & \ddots & \\ 0 & \ddots & \ddots & \ddots & \ddots & \ddots & \ddots & \ddots & \ddots & \vdots \\ \vdots & \ddots & \ddots & \ddots & \ddots & \ddots & \ddots & \ddots & \ddots & \\ 0 & \ddots & \ddots & \ddots & \ddots & \ddots & \cdots & a_{n-3} & A_{n-3} & C_{n-3} \\ 0 & \cdots & \ddots & \ddots & \ddots & \ddots & \cdots & d_{n-2} & a_{n-2} & A_{n-2} \\ A_{n-1} & 0 & \cdots & \ddots & \ddots & \ddots & \cdots & b_{n-1} & d_{n-1} & a_{n-1} \\ a_n & A_n & 0 & \cdots & \cdots & & D_n & B_n & b_n & d_n \end{bmatrix}, \tag{1.2}$$

---

[1] Home address: Mathematics Department, Faculty of Science, Mansoura University, Mansoura, 35516, Egypt. E-mail:abibka@mans.edu.eg

$X = (x_1, x_2, ..., x_n)^T$, $R = (R_1, R_2, ..., R_n)^T$, $n \geq 8$.

In paper [3], Karawia give an efficient symbolic algorithm to obtain the inverse of heptadiagonal matrix of the form (1.2) and then the solution of the systems (1.1) based on LU decomposition. There are many special cases of cyclic heptadiagonal linear systems (See [4-11])

Recently in [12], a new efficient computational algorithm is presented for solving nearly penta-diagonal linear systems based on the use of any penta-diagonal linear solver. In this paper we are going to compute the solution of a general cyclic heptadiagonal system of the form (1.1) without imposing any restrictive conditions on the elements of the matrix H in (1.2). Our approach is mainly based on getting the natural generalization of the algorithm presented in [12]. The development of a symbolic algorithm is considered in order to remove all cases where the numerical algorithm fails.

The paper is organized as follows. In Section 2, symbolic computational algorithm for the solution of heptadiagonal linear system, that will not break, is constructed. In Section 3, Sherman–Morrison–Woodbury formula is given. An illustrative example is given in section 4. In Section 5, Conclusions of the work are presented.

## 2. Heptadiagonal linear solver

In this section we shall focus on the construction of new symbolic computational algorithms for computing the solution of general heptadiagonal linear system of the form:

$$H_h X_h = R_h, \tag{2.1}$$

where

$$H_h = \begin{bmatrix} d_1 & a_1 & A_1 & C_1 & 0 & 0 & \cdots & & 0 & 0 \\ b_2 & d_2 & a_2 & A_2 & C_2 & \ddots & \ddots & \cdots & & 0 \\ B_3 & b_3 & d_3 & a_3 & A_3 & \ddots & \ddots & \ddots & \cdots & \\ D_4 & B_4 & b_4 & d_4 & a_4 & \ddots & \ddots & \ddots & \ddots & \\ 0 & \ddots & \ddots & \ddots & \ddots & \ddots & \ddots & \ddots & \ddots & \vdots \\ \vdots & \ddots & \ddots & \ddots & \ddots & \ddots & \ddots & \ddots & \ddots & \\ 0 & \ddots & \ddots & \ddots & \ddots & \ddots & \cdots & a_{m-3} & A_{m-3} & C_{m-3} \\ 0 & \cdots & \ddots & \ddots & \ddots & \ddots & \cdots & d_{m-2} & a_{m-2} & A_{m-2} \\ 0 & 0 & \cdots & \ddots & \ddots & \ddots & \cdots & b_{m-1} & d_{m-1} & a_{m-1} \\ 0 & 0 & 0 & \cdots & \cdots & & D_m & B_m & b_m & d_m \end{bmatrix} \tag{2.2}$$

$X_h = (x_{h_1}, x_{h_2}, ..., x_{h_m})^T$, $R = (R_{h_1}, R_{h_2}, ..., R_{h_m})^T$, $m \geq 4$.

Firstly we begin with computing the LU factorization of the matrix $H_h$. It is as in the following:

$$H_h = LU \tag{2.3}$$

where

$$
L=\begin{bmatrix}
1 & 0 & 0 & 0 & 0 & \cdots & & 0 & 0 & 0 \\
f_2 & 1 & 0 & 0 & 0 & \cdots & & 0 & 0 & 0 \\
e_3 & f_3 & 1 & 0 & 0 & \cdots & & 0 & 0 & 0 \\
\frac{D_4}{\alpha_1} & e_4 & f_4 & 1 & 0 & \cdots & & 0 & 0 & 0 \\
0 & \ddots & \ddots & \ddots & \ddots & \ddots & \vdots & \vdots & \vdots & \vdots \\
\vdots & \ddots & \ddots & \ddots & \ddots & \ddots & \ddots & \vdots & \vdots & \vdots \\
0 & 0 & 0 & \ddots & \ddots & \ddots & \ddots & \ddots & \vdots & \vdots \\
0 & 0 & 0 & \cdots & \ddots & \ddots & \ddots & \ddots & \ddots & \vdots \\
0 & 0 & 0 & \cdots & & \frac{D_{m-1}}{\alpha_{m-4}} & e_{m-1} & f_{m-1} & 1 & 0 \\
0 & 0 & 0 & \cdots & & 0 & \frac{D_m}{\alpha_{m-3}} & e_m & f_m & 1
\end{bmatrix} \quad (2.4)
$$

and

$$
U=\begin{bmatrix}
\alpha_1 & g_1 & z_1 & C_1 & 0 & \cdots & & 0 & 0 & 0 \\
0 & \alpha_2 & g_2 & z_2 & \ddots & \ddots & & 0 & 0 & 0 \\
0 & 0 & \alpha_3 & g_3 & \ddots & \ddots & & 0 & 0 & 0 \\
\vdots & & \ddots & \ddots & \ddots & \ddots & \ddots & \vdots & \vdots & \vdots \\
\vdots & & & \ddots & \ddots & \ddots & \ddots & \ddots & & 0 \\
\vdots & & & & \ddots & \ddots & \ddots & \ddots & \ddots & \\
\vdots & & & & & \ddots & \dots & \dots & \dots & C_{m-3} \\
\vdots & & & & & & 0 & \alpha_{m-2} & g_{m-2} & z_{n-2} \\
\vdots & & & & & & 0 & 0 & \alpha_{m-1} & g_{m-1} \\
0 & 0 & 0 & \cdots & \cdots & \cdots & & 0 & 0 & \alpha_m
\end{bmatrix} \quad (2.5)
$$

The elements in the matrices L and U in (2.4) and (2.5) satisfy:

$$
\alpha_i=\begin{cases}
d_1 & \text{if } i=1 \\
d_2-f_2g_1 & \text{if } i=2 \\
d_3-e_3z_1-f_3g_2 & \text{if } i=3 \\
d_i-\dfrac{D_i}{\alpha_{i-3}}C_{i-3}-e_iz_{i-2}-f_ig_{i-1} & \text{if } i=4,5,\ldots,m,
\end{cases} \quad (2.6)
$$

$$
f_i=\begin{cases}
\dfrac{b_2}{\alpha_1} & \text{if } i=2 \\
\dfrac{b_3-e_3g_1}{\alpha_2} & \text{if } i=3 \\
\dfrac{1}{\alpha_{i-1}}\left(b_i-\dfrac{D_i}{\alpha_{i-3}}z_{i-3}-e_ig_{i-2}\right) & \text{if } i=4,5,\ldots,m,
\end{cases} \quad (2.7)
$$

$$e_i = \begin{cases} \dfrac{B_3}{\alpha_1} & \text{if } i = 3 \\ \dfrac{1}{\alpha_{i-2}}\left(B_i - \dfrac{D_i}{\alpha_{i-3}} g_{i-3}\right) & \text{if } i = 4,5,\ldots,m, \end{cases} \tag{2.8}$$

$$g_i = \begin{cases} a_1 & \text{if } i = 1 \\ a_2 - f_2 z_1 & \text{if } i = 2 \\ a_i - f_i z_{i-1} - e_i C_{i-2} & \text{if } i = 3,4,\ldots,m-1, \end{cases} \tag{2.9}$$

$$z_i = \begin{cases} A_1 & \text{if } i = 1 \\ A_i - f_i C_{i-1} & \text{if } i = 2,3,\ldots,m-2. \end{cases} \tag{2.10}$$

We also have:

$$\det H_h = \prod_{i=1}^{n} \alpha_i . \tag{2.11}$$

Then, the solution is given by

$$x_{h_i} = \begin{cases} \dfrac{Q_{h_m}}{\alpha_m} & \text{if } i = m \\ \dfrac{1}{\alpha_{m-1}}\left(Q_{h_{m-1}} - g_{m-1} x_{h_m}\right) & \text{if } i = m-1 \\ \dfrac{1}{\alpha_{m-2}}\left(Q_{h_{m-2}} - g_{m-2} x_{h_{m-1}} - z_{m-2} x_{h_m}\right) & \text{if } i = m-2 \\ \dfrac{1}{\alpha_i}\left(Q_{h_i} - g_i x_{h_{i+1}} - z_i x_{h_{i+2}} - C_i x_{h_{i+3}}\right) & \text{if } i = m-3, m-4, \ldots, 1, \end{cases} \tag{2.12}$$

where

$$Q_{h_i} = \begin{cases} R_{h_1} & \text{if } i = 1 \\ R_{h_2} - f_2 Q_{h_1} & \text{if } i = 2 \\ R_{h_3} - e_3 Q_{h_1} - f_3 Q_{h_2} & \text{if } i = 3 \\ R_{h_i} - \dfrac{D_i}{\alpha_{i-3}} Q_{h_{i-3}} - e_i Q_{h_{i-2}} - f_i Q_{h_{i-1}} & \text{if } i = 4,5,\ldots,1. \end{cases} \tag{2.13}$$

At this point it is convenient to formulate our first result. It is a symbolic algorithm for computing the solution of a heptadiagonal linear system of the form (2.1).

**Algorithm 2.1.** To compute the solution of a heptadiagonal linear system of the form (2.1), we may proceed as follows:

**Step 1:** Set $\alpha_1$=$d_1$. If $\alpha_1$=0 then $\alpha_1$= t end if. Set $g_1$=$a_1$, $z_1$=$A_1$, $k_1$=$A_{n-1}$/ $\alpha_1$, $f_2$=$b_2$/ $\alpha_1$, $e_3$=$B_3$/ $\alpha_1$, $\alpha_2$=$d_2$-$f_2$ *$g_1$. If $\alpha_2$=0 then $\alpha_2$= t end if. Set $g_2$=$a_2$-$f_2$*$z_1$, $f_3$=($b_3$-$e_3$*$g_1$)/$\alpha_2$, $\alpha_3$=$d_3$-$e_3$*$z_1$-$f_3$*$g_2$. If $\alpha_3$=0 then $\alpha_3$= t end if.

**Step 2:** Compute and simplify:

For i from 4 to m do

$e_i$=($B_i$-$D_i$* $g_{i-3}$/ $\alpha_{i-3}$)/ $\alpha_{i-2}$

$f_i$=($b_i$- $D_i$* $z_{i-3}$/ $\alpha_{i-3}$-$e_i$*$g_{i-2}$)/ $\alpha_{i-1}$

$z_{i-2}$=$A_{i-2}$-$f_{i-2}$*$C_{i-3}$

$g_{i-1}$=$a_{i-1}$-$f_{i-1}$*$z_{i-2}$-$e_{i-1}$*$C_{i-3}$

$\alpha_i=(d_i- D_i* C_{i-3}/ \alpha_{i-3}-e_i z_{i-2}-f_i*g_{i-1})$

If $\alpha_i=0$ then $\alpha_i= t$ end if

End do

**Step 3:** Compute det $H_h = \left(\prod_{i=1}^{n}\alpha_i\right)_{t=0}$ .

**Step 4:** If det $H_h \neq 0$ then do

**Step 5:** Set $Q_{h1}=R_{h1}$, $Q_{h2}=R_{h2}-f_2*Q_{h1}$, $Q_{h3}=R_{h3}-e_3*Q_{h1}-f_3*Q_{h2}$,

Compute and simplify:

For i from 4 to m do

$Q_{hi}=R_{hi}-D_i*Q_{hi-3}/ \alpha_{i-3}- e_i*Q_{hi-2}-f_i*Q_{hi-1}$

End do

**Step 6:** Set $x_{hm}=Q_{hm}/ \alpha_m$, $x_{hm-1}=(Q_{hm-1}-g_{m-1}*x_{hm})/ \alpha_{m-1}$, $x_{hm-2}=(Q_{hm-2}-g_{m-2} *x_{hm-1}-z_{m-2}*x_{hm})/ \alpha_{m-2}$ ,

Compute and simplify:

For i from m-3 by -1 to 1 do

$x_{hi}=(Q_{hi}-g_i *x_{hi+1}-z_i*x_{hi+2}-C_i*x_{hi+3})/ \alpha_i$ .

End do

**Step 8:** Compute and simplify the solution:

For i from 1 to m do

$X_{hi}=(x_{hi})_{t=0.}$

End do

Else

OUTPUT("The matrix $H_h$ is singular"); Stop.

End If

The new algorithm 2.1 is very useful to check the nonsingularity of the matrix $H_h$.

## 3. Sherman–Morrison–Woodbury formula[13]

In this section, we are going to formulate a new computational algorithm for solving cyclic heptadiagonal linear systems of the form (1.1) based on the previous heptadiagonal linear solver. The heptadiagonal linear system of the form (1.1) can be written in the form:

$$\begin{pmatrix} M_1 & V \\ U^T & M_2 \end{pmatrix}\begin{pmatrix} x' \\ x'' \end{pmatrix}=\begin{pmatrix} R' \\ R'' \end{pmatrix} \tag{3.1}$$

where

$$M_1=\begin{bmatrix} d_1 & a_1 & A_1 & C_1 & 0 & 0 & \cdots & & 0 & 0 \\ b_2 & d_2 & a_2 & A_2 & C_2 & \ddots & \ddots & \cdots & & 0 \\ B_3 & b_3 & d_3 & a_3 & A_3 & \ddots & \ddots & \ddots & \cdots & \\ D_4 & B_4 & b_4 & d_4 & a_4 & \ddots & \ddots & \ddots & \ddots & \\ 0 & \ddots & \ddots & \ddots & \ddots & \ddots & \ddots & \ddots & \ddots & \vdots \\ \vdots & \ddots & \ddots & \ddots & \ddots & \ddots & \ddots & \ddots & \ddots & \\ 0 & \ddots & \ddots & \ddots & \ddots & \ddots & \ddots & a_{n-5} & A_{n-5} & C_{n-5} \\ 0 & \cdots & \ddots & \ddots & \ddots & \ddots & \ddots & d_{n-4} & a_{n-4} & A_{n-4} \\ 0 & 0 & \cdots & \ddots & \ddots & \ddots & \ddots & b_{n-3} & d_{n-3} & a_{n-3} \\ 0 & 0 & 0 & \cdots & \cdots & & D_{n-2} & B_{n-2} & b_{n-2} & d_{n-2} \end{bmatrix}, M_2=\begin{bmatrix} d_{n-1} & a_{n-1} \\ b_n & d_n \end{bmatrix},$$

$$V=\begin{bmatrix} v_1^T \\ v_2^T \end{bmatrix}^T=\begin{bmatrix} B_1 & 0 & 0 & \cdots & \cdots & \cdots & C_{n-4} & A_{n-3} & a_{n-2} \\ b_1 & B_2 & 0 & \cdots & \cdots & \cdots & 0 & C_{n-3} & A_{n-2} \end{bmatrix}^T,$$

$$U^T=\begin{bmatrix} A_{n-1} & 0 & 0 & \cdots & \cdots & \cdots & D_{n-1} & B_{n-1} & b_{n-1} \\ a_n & A_n & 0 & \cdots & \cdots & \cdots & 0 & D_n & B_n \end{bmatrix},$$

$x'=(x_1,x_2,...,x_{n-2})^T$, $x''=(x_{n-1},x_n)^T$, $R'=(R_1,R_2,...,R_{n-2})^T$, and $R''=(R_{n-1},R_n)^T$.

Thus (3.1) is equivalent to

$$\begin{aligned} M_1x'+V\ x''&=R' \\ U^Tx'+M_2x''&=R'' \end{aligned} \tag{3.2}$$

Assume that $M_2$ is nonsingular. After elimination of $x''$ from (3.2), we get the linear systems

$$M\ x'=\hat{R} \tag{3.3}$$

where $M=M_1-VM_2^{-1}U^T$, $\hat{R}=R'-VM_2^{-1}R''$.

If we applying the Sherman–Morrison–Woodbury formula to M, we obtain

$$M^{-1}=M_1^{-1}+M_1^{-1}V(M_2-U^TM_1^{-1}V)^{-1}U^TM_1^{-1} \text{ and } x'=M^{-1}\hat{R}=y+M_1^{-1}V(M_2-U^TM_1^{-1}V)^{-1}U^Ty$$

where y is the solution of $M_1y=\hat{R}$. It is clear that the solution $x''$ can be found from the above formula by successive calculation of the expressions

$$y=M_1^{-1}\hat{R},\ M_1^{-1}V,\ U^TM_1^{-1}V,\ (M_2-U^TM_1^{-1}V)^{-1}, \text{and } (M_2-U^TM_1^{-1}V)^{-1}U^Ty.$$

The main part of above calculations is in finding the first two expressions, which is equivalent to solving three (n _ 2)-by-(n _ 2) heptadiagonal linear systems with the same coefficient matrix $M_1$ and different right-hand sides. After finding of $x'$, we can get $x''$ from the second equation of (3.2) by formula

$$x''=M_2^{-1}(R''-U^Tx')$$

At this point it is convenient to formulate our second result. It is a symbolic algorithm for computing the solution of a cyclic heptadiagonal linear system of the form (1.1) and can be considered as natural generalization of the symbolic algorithm 1 in [12].

**Algorithm 3.1.** To compute the solution of a cyclic heptadiagonal linear system of the form (1.1), we may proceed as follows:

**Step 1:** Find $M_1$, $M_2$, $U^T$, V, $R', R''$, and $\hat{R}=R'-VM_2^{-1}R''$.

**Step 2:** Solve $M_1y=\hat{R}, M_1q_1=v_1$, and $M_1q_2=v_2$ by algorithm 2.1, then obtain y and $M_1^{-1}v=(q_1,q_2)$.

**Step 3:** Compute $x'=y+(q_1,q_2)(M_2-U^T(q_1,q_2))^{-1}U^Ty$, $x''=M_2^{-1}(R''-U^Tx')$.

**Step 4:** Compute the solution $x=\begin{pmatrix} x' \\ x'' \end{pmatrix}_{t=0}$.

Three systems $M_1y=\hat{R}, M_1q_1=v_1$, and $M_1q_2=v_2$ in algorithm 3.1 can be solved in parallel.

## 4. An illustrative example

In this section, we give simple numerical experiment to illustrate the effectiveness of our algorithms.

**Example 4.1**. Consider the 10-by-10 cyclic heptadiagonal systems coming from [3]

$$\begin{bmatrix} 1 & -1 & 1 & -2 & 0 & 0 & 0 & 0 & 2 & -1 \\ 1 & 1 & 1 & 1 & -1 & 0 & 0 & 0 & 0 & 1 \\ 2 & 1 & -1 & 1 & 2 & 3 & 0 & 0 & 0 & 0 \\ 2 & -2 & 3 & 1 & 5 & -6 & 0 & 0 & 0 & 0 \\ 0 & 1 & 1 & 1 & 1 & 1 & 1 & 2 & 0 & 0 \\ 0 & 0 & -1 & -1 & -1 & -1 & -1 & -1 & 1 & 0 \\ 0 & 0 & 0 & 2 & 2 & 2 & 2 & 3 & 1 & -3 \\ 0 & 0 & 0 & 0 & -2 & -2 & 1 & 1 & 3 & 5 \\ 3 & 0 & 0 & 0 & 0 & 3 & 1 & 3 & 4 & -1 \\ 2 & 4 & 0 & 0 & 0 & 0 & 2 & 3 & 4 & 1 \end{bmatrix} \begin{bmatrix} x_1 \\ x_2 \\ x_3 \\ x_4 \\ x_5 \\ x_6 \\ x_7 \\ x_8 \\ x_9 \\ x_{10} \end{bmatrix} = \begin{bmatrix} 2 \\ 15 \\ 33 \\ 0 \\ 43 \\ -24 \\ 47 \\ 70 \\ 78 \\ 94 \end{bmatrix} \tag{4.1}$$

**Solution:** The application of Algorithm 3.1 gives

**Step1:** $M_1 = \begin{bmatrix} 1 & -1 & 1 & -2 & 0 & 0 & 0 & 0 \\ 1 & 1 & 1 & 1 & -1 & 0 & 0 & 0 \\ 2 & 1 & -1 & 1 & 2 & 3 & 0 & 0 \\ 2 & -2 & 3 & 1 & 5 & -6 & 0 & 0 \\ 0 & 1 & 1 & 1 & 1 & 1 & 1 & 2 \\ 0 & 0 & -1 & -1 & -1 & -1 & -1 & -1 \\ 0 & 0 & 0 & 2 & 2 & 2 & 2 & 3 \\ 0 & 0 & 0 & 0 & -2 & -2 & 1 & 1 \end{bmatrix}, M_2 = \begin{bmatrix} 4 & -1 \\ 4 & 1 \end{bmatrix},$

$$U^T = \begin{bmatrix} 3 & 0 & 0 & 0 & 0 & 3 & 1 & 3 \\ 2 & 4 & 0 & 0 & 0 & 0 & 2 & 3 \end{bmatrix}, V = \begin{bmatrix} v_1^T \\ v_2^T \end{bmatrix}^T = \begin{bmatrix} 2 & 0 & 0 & 0 & 0 & 1 & 1 & 3 \\ -1 & 1 & 0 & 0 & 0 & 0 & -3 & 5 \end{bmatrix}^T,$$

$$R' = \begin{bmatrix} 2 & 15 & 33 & 0 & 43 & -24 & 47 & 70 \end{bmatrix}^T, R'' = \begin{bmatrix} 78 & 94 \end{bmatrix}^T,$$

$$\hat{R} = R' - VM_2^{-1}R'' = (-33, 7, 33, 0, 43, \frac{-91}{2}, \frac{99}{2}, \frac{-69}{2})^T.$$

**Step2:** $y = \left( \frac{-2814}{199}, \frac{3345}{199}, \frac{2208}{199}, \frac{1308}{199}, \frac{2654}{199}, \frac{4442}{597}, \frac{15739}{597}, \frac{-7685}{398} \right)^T,$

$$M_1^{-1}V = (q_1, q_2) = \begin{pmatrix} \frac{242}{199} & \frac{-98}{199} & \frac{-150}{199} & \frac{-104}{199} & \frac{-110}{199} & \frac{-212}{597} & \frac{-184}{597} & \frac{297}{199} \\ \frac{6}{199} & \frac{861}{199} & \frac{-132}{199} & \frac{-394}{199} & \frac{142}{199} & \frac{-895}{597} & \frac{4630}{597} & \frac{-861}{199} \end{pmatrix}^T.$$

**Step3:** $x' = y + (q_1, q_2)(M_2 - U^T(q_1, q_2))^{-1}U^T y = (1, 2, 3, 4, 5, 6, 7, 8),$

$$x'' = M_2^{-1}(R'' - U^T x') = (9,\ 10).$$

**Step 4:** The solution of (4.1) is

$$x = \begin{pmatrix} x' \\ x'' \end{pmatrix} = (1,2,3,4,5,6,7,8,9,10)^T.$$

## 5. Conclusions

In this paper, we derived a computational algorithm for solving the cyclic heptadiagonal linear systems. Since the algorithm uses the symbolic algorithm 2.1 method for the LU factorization of the heptadiagonal matrix, the factorization never suffers from breakdown, and this leads to the fast and reliable solution of cyclic heptadiagonal linear systems. The realization of the method needs O(n) operations. The algorithms are natural generalizations of some algorithms in current use.